\documentclass{article}
\usepackage{amsmath}
\usepackage{graphicx}
\newtheorem{theorem}{Theorem}

\begin{document}
\title{A new kernel estimator of hazard ratio and its asymptotic mean squared error}
\author{
Taku MORIYAMA\\
Graduate School of Mathematics\\
Yoshihiko MAESONO\\
Faculty of Mathematics\\\vspace{5mm}\\
Kyushu University, Motooka 744, Fukuoka 819-0395, Japan}
\date{}

\maketitle

\begin{abstract}
The hazard function is a ratio of a density and survival function, and it is a basic tool of the survival analysis. In this paper we propose a kernel estimator of the hazard ratio function, which are based on a modification of \'{C}wik and Mielniczuk \cite{cwik1989estimating}'s method. We study nonparametric estimators of the hazard function and compare those estimators by means of asymptotic mean squared error ($AMSE$).  We obtain asymptotic bias and variance of the new estimator, and compare them with a naive estimator.  The asymptotic variance of the new estimator is always smaller than the naive estimator's, so we also discuss an improvement of $AMSE$ using Terrell and Scott \cite{terrell1980improving}'s bias reduction method. The new modified estimator ensures the non-negativity, and we demonstrate the numerical improvement.
\end{abstract}

\section{Introduction}

Rosenblatt \cite{rosenblatt1956remarks} has proposed a kernel estimator of a probability density function $f(\cdot)$.  After that, many researchers have discussed various kernel estimators for distribution, regression, hazard functions etc.  Since most of the kernel type estimators are biased, many papers have studied bias reductions.  Those are based on a higher order kernel, transformations of estimators etc. Terrell and Scott \cite{terrell1980improving} have proposed efficient transformation.  We will apply their method to the kernel estimator in Section 4.  On the other hand, \'{C}wik and Mielniczuk \cite{cwik1989estimating} have discussed a new kernel estimator of a density ratio which they called `direct' estimator.  Asymptotic mean squared error ($AMSE$) of the direct estimator is different from a naive estimator, and the variance is always smaller.  Here let us introduce their direct estimator.

Let $X_1,X_2,\cdots,X_n$ be independently and identically distributed ($i.i.d.$) random variables with a distribution function $F(\cdot)$, and $Y_1, Y_2 ,\cdots, Y_n$ be $i.i.d.$ random variables with a distribution function $G(\cdot)$.  $f(\cdot)$ and $g(\cdot)$ are the density functions of $F(\cdot)$ and $G(\cdot)$, and we assume $g(x_0) \neq 0 ~(x_0 \in {\bf R})$.  A naive estimator of $f(x_0)/g(x_0)$ is defined as replacing $f(x_0)$ and $g(x_0)$ by the kernel density estimators $\widehat{f}(x_0)$ and $\widehat{g}(x_0)$ where 
$$\widehat{f}(x_0)=\frac{1}{nh} \sum_{i=1}^n K\left(\frac{x_0 -X_i}{h} \right)\hspace{4mm}{\rm and}\hspace{4mm} \widehat{g}(x_0)=\frac{1}{nh} \sum_{i=1}^n K\left(\frac{x_0 -Y_i}{h} \right).$$
$K(\cdot)$ is a kernel function and $h$ is a bandwidth which satisfies $h \rightarrow 0$ and $n h \rightarrow \infty ~(n \rightarrow \infty)$.  We call $\widehat{f}(x_0)/\widehat{g}(x_0)$ `indirect' estimator.  The direct type estimator has been proposed by \'{C}wik and Mielniczuk \cite{cwik1989estimating}, which is given by
$$\widehat{\frac{f}{g}}(x_0)=\frac{1}{h} \int_{-\infty}^{\infty} K\left( \frac{G_{n}(x_0) -G_{n}(y) }{h} \right) dF_n(y),$$
where $F_n(\cdot)$ and $G_{n}(\cdot)$ are the empirical distribution functions of $X_1,\cdots,X_n$ and $Y_1,\cdots, Y_n$, respectively. Chen, Hsu, and Liaw \cite{chen2009kernel} have obtained an explicit form of its $AMSE$.

Extending their method, we propose a new `direct' estimator of the hazard ratio function in Section 2 and investigate its $AMSE$.  We will compare the naive and direct kernel estimators in Section 3. It is demonstrated that our direct estimator performs asymptotically better, especially when the underlying distribution is exponential or gamma.  Though the bias term of the direct estimator is large in some cases, the asymptotic variance is always smaller.  As mentioned before, there are many bias reduction methods and then we discuss the bias reduction in Section 4.  Applying the Terrell and Scott \cite{terrell1980improving}'s method, we will show that the modified direct estimator gives us good performance both theoretically and numerically.  Proofs of them are given in Appendices.

\section{Hazard ratio estimators and asymptotic properties}
\label{sect:Hazard ratio estimators and asymptotic properties}

\subsection{Hazard ratio estimators} 

Let us assume that the density function $f(\cdot)$ of $X_i$ satisfies $f(x_0) \neq 0 ~(x_0 \in {\bf R})$.  The hazard ratio function is defined as
$$H(x) = \frac{f(x)}{1-F(x)}.$$
The meaning of $H(x) dx$ is conditional probability of `death' in $[x, x+dx]$ given a survival to $x$. This is a fundamental measure in survival analysis, and the range of applications is wide.  For example, actuaries call it ``force of mortality'' and use it to estimate insurance money.  In extreme value theory, it is called ``intensity function'' and used to determine the form of an extreme value distribution (see Gumbel \cite{gumbel1958statistics})

The naive nonparametric estimator of $H(x_0)$ is given by Watson \cite{watson1964hazard}
$$\widetilde{H}(x_0)=\frac{\widehat{f}(x_0)}{1-\widehat{F}(x_0)},$$
where 
$$\widehat{f}(x_0)=\frac{1}{nh} \sum_{i=1}^n K\left(\frac{x_0 -X_i}{h} \right)\hspace{5mm}{\rm and}\hspace{5mm}\widehat{F}(x_0)=\frac{1}{n} \sum_{i=1}^n W\left(\frac{x_0 -X_i}{h} \right).$$
Here $K(\cdot)$ is the kernel function and $W(\cdot)$ is an integral of $K(\cdot)$
$$W(u) = \int _{-\infty}^u K(t) dt.$$
Using properties of the kernel density estimator, Murthy \cite{murthy1965estimation} proves the consistency and asymptotic normality of $\widetilde{H}(x_0)$. In random censorship model, Tanner, Wong et al. \cite{tanner1983estimation} proves them, using the H\'{a}jek's projection method. Patil \cite{patil1993least} gives its mean integrated squared error ($MISE$) and discusses the optimal bandwidth in both uncensored and censored settings.  For dependent data, $MSE$ of the indirect estimator is obtained by Quintela-del R\'{i}o \cite{quintela2007plug}.  Using Vieu \cite{vieu1991quadratic}'s results, he has obtained modified $MISE$ which avoids that the denominator equal to $0$.  In this paper, we assume that the support of the kernel $K(\cdot)$ is $[-d, d] ~(d>0)$ and the data do not have a censoring.

Extending the idea of \'{C}wik and Mielniczuk \cite{cwik1989estimating}, we propose a new `direct' estimator of the hazard ratio function as follows:
$$\widehat{H}(x_0) =\frac{1}{h} \int_{-\infty}^{\infty} K\left( \frac{y- \eta_n(y) - \left( x_0 -\eta_n(x_0) \right) }{h} \right) dF_n(y),$$
where
$$\eta_n (y)= \int_{-\infty}^y F_n(u) du = n^{-1}\sum I(X_i \le y) (y- X_i).$$  It is easy to see that $\widehat{H}(x)$ is a smooth function.  We will discuss its asymptotic properties.

\subsection{Asymptotic properties} 

For brevity, we use a notation
$$A_{i,j} = \int_{-\infty}^{\infty} K^i(u) u^j du.$$
The proofs of the following theorems are given in Appendix.  For the direct hazard estimator, we have the following $AMSE$.

\begin{theorem} Let us assume that ${\rm (i)}$ $f(\cdot)$ is 4 times continuously differentiable at $x_0$, ${\rm(ii)}$ $K$ is symmetric and the support is $[-d, d]$ for some positive number $d>0$ and ${\rm(iii)}$ $A_{1,4}$ and $A_{2,0}$ are bounded. Then, the $MSE$ of $\widehat{H}(x_0)$ is given by
\begin{eqnarray}
&&E\left\{ \widehat{H}(x_0) -\left[\frac{f}{1-F}\right](x_0) \right\}^2\nonumber\\
&=&\frac{h^4}{4} A_{1,2}^2 \left[\frac{\{(1 -F)\{(1-F)f'' +4f f' \} +3f^3\}^2}{(1 -F)^{10}}\right](x_0)+\frac{A_{2,0}}{n h}\left[\frac{f}{1-F}\right] (x_0)  \label{variance1}\\
&&+O\left(h^6+\frac{1}{n h^{1/2}}\right).\nonumber
\end{eqnarray}
\end{theorem}

On the other hand, under some regularity conditions, $MSE$ of $\widetilde{H}(x_0)$ is given by Patil \cite{patil1993least} as follows:
\begin{eqnarray}
&&E\left\{ \widetilde{H}(x_0) -\left[\frac{f}{1-F}\right](x_0) \right\}^2\nonumber\\
&=&\frac{h^4}{4} A_{1,2}^2 \left[\frac{\{(1-F)f'' +f f' \}^2}{(1 -F)^{4}}\right](x_0) +\frac{A_{2,0}}{n h}\left[\frac{f}{(1-F)^2}\right] (x_0) \label{variance2}\\
&&+O\left(h^6+\frac{1}{n h^{1/2}}\right).\nonumber
\end{eqnarray}
The asymptotic variances are second terms of (\ref{variance1}) and (\ref{variance2}), and then the direct estimator has small variance because of $0<1-F(x_0)<1$.  Taking the derivative of the $AMSE$ with respect to $h$, we have an optimal bandwidth $h =h^*$ of $\widehat{H}(x_0)$, where
$$h^* = n^{-1/5} \left( \frac{A_{2,0}}{A_{1,2}^2}\left[\frac{(1 -F)^9 f}{\{(1 -F)\{(1-F)f'' +4f f' \} +3f^3\}^2}\right](x_0) \right)^{1/5}.$$
Though $h^*$ depends on unknown functions, it is available by replacing them by their estimators (plug-in method).  Similarly, the following optimal bandwidth of the indirect $\widetilde{H}(x_0)$ is obtained
$$h^{**} = n^{-1/5} \left( \frac{A_{2,0}}{A_{1,2}^2}\left[\frac{(1 -F)^2 f}{\{(1-F)f'' +f f' \}^2}\right](x_0) \right)^{1/5}.$$

Further we can show the asymptotic normality of the direct $\widehat{H}$.
\begin{theorem} Let us assume the same conditions of Theorem 1 and $h =cn^{-d} (0<c, ~\frac15\leq d<\frac12$. Then, the following asymptotic normality of $\widehat{H}(x_0)$ holds.
$$
\sqrt{n h} \left\{ \widehat{H}(x_0) -\left[\frac{f}{1-F}\right](x_0)-h^2B_1 \right\} \xrightarrow{d} N(0,V_1),
$$
where 
$$
B_1 = \frac{A_{1,2}}{2} \left[\frac{(1 -F)\{(1-F)f'' +4f f' \} +3f^3}{(1 -F)^{5}}\right](x_0)\hspace{3mm}{\rm and}\hspace{3mm} ~V_1 = A_{2,0} \left[\frac{f}{1-F}\right] (x_0).
$$
\end{theorem}
Using the Slutsky's theorem, the asymptotic normality of the indirect estimator is easily obtained.

The direct estimator is superior in the sense of the asymptotic variance, so we will consider the bias reduction in Section 4. We have the following higher-order asymptotic bias.
\begin{theorem} Let us assume that ${\rm(i')}$ $f$ is 6 times continuously differentiable at $x_0$, ${\rm(ii')}$ $K$ is symmetric and the support is $[-d, d]$ for some positive number $d>0$ and ${\rm(iii')}$ $A_{1,6}$ is bounded. Then, the higher-order asymptotic bias of $\widehat{H}(x_0)$ is 
$$E\left\{ \widehat{H}(x_0) -\left[\frac{f}{1-F}\right](x_0) \right\} = h^2a_2(x_0)+h^4a_4(x_0) +O\left(h^6 +n^{-1}\right),$$
where
\begin{eqnarray*}
a_2(x_0)&=&\frac{A_{1,2}}{2}\left[\frac{-m^2 m''' +4m m' m'' -3(m')^3}{m^5}\right] (x_0),\\
a_4(x_0)&=&\frac{A_{1,4}}{24} \biggl[\frac{-60 m^2 (m')^2 m''' +15 m^3 m'' m''' +11 m^3 m' m^{(4)} -m^4 m^{(5)} }{m^9}\\
&&\hspace*{20mm}+ \frac{210m (m')^3 m'' -73 m^2 m' (m'')^2 -105(m')^5}{m^9} \biggl](x_0)
\end{eqnarray*}
and $m(x)=1-F(x)$.
\end{theorem}

\section{Comparison of kernel hazard estimators}
\label{sect:Comparison of kernel hazard estimators}

In some special cases, we will investigate the $AMSE$ values of the direct $\widehat{H}(x_0)$ and indirect $\widetilde{H}(x_0)$. We show that the proposed estimator $\widehat{H}(x_0)$ performs asymptotically better when $F(\cdot)$ is the exponential or gamma distribution. 

Here we consider that $F(\cdot)$ is the exponential, uniform, gamma, weibull or beta distribution. The cumulative distribution function of the exponential distribution $Exp(1/\lambda)$ is $F(x) =1- \exp(-\lambda x)$ and the hazard ratio is constant, that is $H(x) = \lambda$.  It is one of the most common and an important model of the survival analysis.  When $F(\cdot)$ is the exponential, the asymptotic biases of $\widetilde{H}(x_0)$ and $\widehat{H}(x_0)$ vanish and the $AMSE$s are
$$AMSE \left[ \widehat{H}(x_0) \right] = \frac{\lambda}{n h} A_{2,0}\hspace{2mm}<\hspace{2mm} AMSE \left[ \widetilde{H}(x_0) \right] = \frac{\lambda}{n h} \exp(\lambda x_0) A_{2,0}. $$
Thus the proposed estimator is always asymptotically better regardless of the parameter $\lambda$ and the point $x_0$.

Next, let us assume that $F(\cdot)$ is the uniform distribution ($F(x) = x/b ~(0<x<b)$), then the hazard ratio is $H(x) = (b -x)^{-1}$. In this model, the hazard ratio increases drastically in the tail area.  These $AMSE$s are the followings
\begin{eqnarray*}
AMSE \left[ \widehat{H}(x_0) \right]& =& \frac{h^4}{4} A_{1,2}^2 \frac{9 b^4}{(b- x_0)^{10}} +\frac{1}{n h}\frac{1}{b -x_0} A_{2,0} \\ 
AMSE \left[ \widetilde{H}(x_0) \right] &=& \frac{1}{n h} \frac{b}{(b- x_0)^2} A_{2,0}.
\end{eqnarray*}
We find that the asymptotic bias of $\widetilde{H}(x_0)$ vanishes and the variance of $\widehat{H}(x_0)$ is smaller. Their asymptotic performance depends on $x_0$ and $b$, but $AMSE$ of the proposed $\widehat{H}(x_0)$ is smaller when the life span $b$ is large.

Last, we consider that $F$ is the gamma $\Gamma(p,100)$, weibull $W(q,100)$ or beta distribution $100 \times B(r,s)$ where $p$, $q$, $r$ and $s$ are their shape parameters. Their scale $\sigma=100$ is moderate as life span. $\Gamma(p,\sigma)$ is the distribution of a sum of $p$ $i.i.d.$ random variables of $Exp(\sigma)$, and so it is most important case.  Its asymptotic squared bias, variance and $AMSE$ values  for some fixed points $x_0$ are found in Table 1, where we omit power terms of $h$.  $\widehat{\gamma}_{\varepsilon}^{[p]}$ and $\widetilde{\gamma}_{\varepsilon}^{[p]}$ represent those values of $\widehat{H}(x_0)$ and $\widetilde{H}(x_0)$ when $x_0$ is $\varepsilon$-th quantile of $\Gamma(p,100)$.  The kernel is the Epanechnikov with $A_{1,2}=1/5, A_{2,0}=3/10$ and $h = n^{-1/5}$.  The coefficients $n^{-4/5}$ are omitted here.

The weibull distribution $W(q,\sigma)$ is also important in survival analysis because the hazard ratio is proportional to polynomial degree $(q-1)$, that is $H(x) = q \sigma^q x^{q-1}$.  $W(1,\sigma)$ is the exponential distribution.  The beta distribution is often used to describe distribution whose support is finite, and it has plentiful shapes.  Table 2, 3 and 4 give each of the least $AMSE$ values using $h^*$ or $h^{**}$ (in Section 2.2) where $w_{\varepsilon}^{[q]}$ and $\beta_{\varepsilon}^{[r,s]}$ stand $AMSE$ values of $\widehat{H}(x_0)$ and $\widetilde{H}(x_0)$ when $x_0$ is $\varepsilon$-th quantile of $W(q,100)$ or $B(r,s)\times100$.   We demonstrate the proposed estimator $\widehat{H}$ performs asymptotically better in most cases of the gamma $\Gamma(p,100)$.  In the weibull cases, the asymptotic performance of our estimator is good and comparable in the beta cases.

\begin{table}
\caption{Asymptotic $Bias^2$, $Var$ and $AMSE$ when $F$ is gamma}
{\resizebox{\textwidth}{!}{\begin{tabular}{|c|c|c|c|c|c|c|c|c|}
  \hline 
   & $Bias^2$ & $Var$ & $AMSE$ & & $Bias^2$ & $Var$ & $AMSE$ \\
  \hline
  $\widehat{\gamma}_{0.05}^{[1/2]}$ & 7.22$\times10^{-2}$ & 4.01$\times10^{-2}$ & 0.112 & $\widehat{\gamma}_{0.1}^{[1/2]}$ & 8.24$\times10^{-5}$ & 2.11$\times10^{-2}$ & 2.11$\times10^{-2}$ \\
  
  $\widetilde{\gamma}_{0.05}^{[1/2]}$ & 6.53$\times10^{-2}$ & 4.22$\times10^{-2}$ & 0.107 & $\widetilde{\gamma}_{0.1}^{[1/2]}$ & 6.72$\times10^{-5}$ & 2.33$\times10^{-2}$ & 2.34$\times10^{-2}$ \\
  \hline
  
  $\widehat{\gamma}_{0.25}^{[1/2]}$ & 1.33$\times10^{-8}$ & 9.52$\times10^{-3}$ & 9.52$\times10^{-3}$ & $\widehat{\gamma}_{0.5}^{[1/2]}$ & 2.36$\times10^{-11}$ & 5.65$\times10^{-3}$ & 5.65$\times10^{-3}$\\
  
  $\widetilde{\gamma}_{0.25}^{[1/2]}$ & 7.74$\times10^{-9}$ & 1.27$\times10^{-2}$ & 1.27$\times10^{-2}$ & $\widetilde{\gamma}_{0.5}^{[1/2]}$ & 6.83$\times10^{-12}$ & 1.13$\times10^{-2}$ & 1.13$\times10^{-2}$\\
  \hline
  
  $\widehat{\gamma}_{0.75}^{[1/2]}$ & 5.35$\times10^{-13}$ & 4.29$\times10^{-3}$ & 4.29$\times10^{-3}$ & $\widehat{\gamma}_{0.9}^{[1/2]}$ & 1.69$\times10^{-15}$ & 3.76$\times10^{-3}$ & 3.76$\times10^{-3}$ \\
  
  $\widetilde{\gamma}_{0.75}^{[1/2]}$ & 6.00$\times10^{-14}$ & 1.72$\times10^{-2}$ & 1.72$\times10^{-2}$ & $\widetilde{\gamma}_{0.9}^{[1/2]}$ & 2.93$\times10^{-15}$ & 3.76$\times10^{-2}$ & 3.76$\times10^{-2}$\\
  \hline
  
  $\widehat{\gamma}_{0.95}^{[1/2]}$ & 1.01$\times10^{-12}$ & 3.58$\times10^{-3}$ & 3.58$\times10^{-3}$ & $\widehat{\gamma}_{0.975}^{[1/2]}$ & 1.47$\times10^{-11}$ & 3.46$\times10^{-3}$ & 3.46$\times10^{-3}$ \\
  
  $\widetilde{\gamma}_{0.95}^{[1/2]}$ & 6.93$\times10^{-16}$ & 7.16$\times10^{-2}$ & 7.16$\times10^{-2}$ & $\widetilde{\gamma}_{0.975}^{[1/2]}$ & 2.33$\times10^{-16}$ & 0.139 & 0.139 \\
  
  \hline  \hline
  $\widehat{\gamma}_{0.05}^{[10]}$ & 2.49$\times10^{-18}$ & 1.56$\times10^{-4}$ & 1.56$\times10^{-4}$ & $\widehat{\gamma}_{0.1}^{[10]}$ & 2.16$\times10^{-18}$ & 2.55$\times10^{-4}$ & 2.55$\times10^{-4}$ \\
  
  $\widetilde{\gamma}_{0.05}^{[10]}$ & 7.16$\times10^{-19}$ & 1.64$\times10^{-4}$ & 1.64$\times10^{-4}$ & $\widetilde{\gamma}_{0.1}^{[10]}$ & 1.73$\times10^{-21}$ & 2.83$\times10^{-4}$ & 2.83$\times10^{-4}$ \\
  \hline
  
  $\widehat{\gamma}_{0.25}^{[10]}$ & 2.61$\times10^{-18}$ & 4.77$\times10^{-4}$ & 4.77$\times10^{-4}$ & $\widehat{\gamma}_{0.5}^{[10]}$ & 1.37$\times10^{-17}$ & 7.72$\times10^{-4}$ & 7.72$\times10^{-4}$\\
  
  $\widetilde{\gamma}_{0.25}^{[10]}$ & 2.40$\times10^{-18}$ & 6.36$\times10^{-4}$ & 6.36$\times10^{-4}$ & $\widetilde{\gamma}_{0.5}^{[10]}$ & 7.91$\times10^{-18}$ & 1.54$\times10^{-3}$ & 1.54$\times10^{-3}$\\
  \hline
  
  $\widehat{\gamma}_{0.75}^{[10]}$ & 2.84$\times10^{-16}$ & 1.07$\times10^{-3}$ & 1.07$\times10^{-3}$ & $\widehat{\gamma}_{0.9}^{[10]}$ & 1.19$\times10^{-14}$ & 1.32$\times10^{-3}$ & 1.32$\times10^{-3}$\\
  
  $\widetilde{\gamma}_{0.75}^{[10]}$ & 1.05$\times10^{-17}$ & 4.28$\times10^{-3}$ & 4.28$\times10^{-3}$ & $\widetilde{\gamma}_{0.9}^{[10]}$ & 9.83$\times10^{-18}$ & 1.32$\times10^{-2}$ & 1.32$\times10^{-2}$\\
  \hline
  
  $\widehat{\gamma}_{0.95}^{[10]}$ & 1.84$\times10^{-13}$ & 1.45$\times10^{-3}$ & 1.45$\times10^{-3}$ & $\widehat{\gamma}_{0.975}^{[10]}$ & 2.75$\times10^{-12}$ & 1.56$\times10^{-3}$ & 1.56$\times10^{-3}$ \\
  
  $\widetilde{\gamma}_{0.95}^{[10]}$ & 8.70$\times10^{-18}$ & 2.90$\times10^{-2}$ & 2.90$\times10^{-2}$ & $\widetilde{\gamma}_{0.975}^{[10]}$ & 7.58$\times10^{-18}$ & 6.24$\times10^{-2}$ & 6.24$\times10^{-2}$ \\
  \hline
\end{tabular}}}
\label{tab:table1}
\end{table}

\begin{table}
\caption{$AMSE$ values when $F$ is gamma and $h = h^*$ or $h^{**}$}
{{\resizebox{\textwidth}{!}{\begin{tabular}{|c|c|c|c|c|c|c|c|}
  \hline 
    &   $AMSE$  &  &  $AMSE$ &  &   $AMSE$  &  &  $AMSE$\\
  \hline
  $\widehat{\gamma}_{0.05}^{[1/2]}$   &  7.44$\times10^{-2}$ & $\widehat{\gamma}_{0.1}^{[1/2]}$    &  1.14$\times10^{-2}$ & $\widehat{\gamma}_{0.25}^{[1/2]}$   &   1.06$\times10^{-3}$ & $\widehat{\gamma}_{0.5}^{[1/2]}$   &  1.97$\times10^{-4}$ \\
  
  $\widetilde{\gamma}_{0.05}^{[1/2]}$    &  7.60$\times10^{-2}$  & $\widetilde{\gamma}_{0.1}^{[1/2]}$   &  1.19$\times10^{-2}$   &  $\widetilde{\gamma}_{0.25}^{[1/2]}$  &  1.20$\times10^{-3}$  &  $\widetilde{\gamma}_{0.5}^{[1/2]}$  & 2.67$\times10^{-4}$ \\
  
    \hline
  $\widehat{\gamma}_{0.75}^{[1/2]}$   &   7.40$\times10^{-5}$ &  $\widehat{\gamma}_{0.9}^{[1/2]}$   &   2.11$\times10^{-5}$ & $\widehat{\gamma}_{0.95}^{[1/2]}$   &   7.26$\times10^{-5}$ &  $\widehat{\gamma}_{0.975}^{[1/2]}$  &  1.21$\times10^{-4}$  \\
  
  $\widetilde{\gamma}_{0.75}^{[1/2]}$  &  1.45$\times10^{-4}$  &  $\widetilde{\gamma}_{0.9}^{[1/2]}$  & 1.48$\times10^{-4}$   &  $\widetilde{\gamma}_{0.95}^{[1/2]}$  &   1.86$\times10^{-4}$ & $\widetilde{\gamma}_{0.975}^{[1/2]}$   &   2.54$\times10^{-4}$ \\
  
  \hline\hline
  $\widehat{\gamma}_{0.05}^{[10]}$   &  4.48$\times10^{-7}$ & $\widehat{\gamma}_{0.1}^{[10]}$  &  6.45$\times10^{-7}$ & $\widehat{\gamma}_{0.25}^{[10]}$   &   1.11$\times10^{-6}$ & $\widehat{\gamma}_{0.5}^{[10]}$   &  2.26$\times10^{-6}$ \\
  
  $\widetilde{\gamma}_{0.05}^{[10]}$   &  3.64$\times10^{-7}$  &  $\widetilde{\gamma}_{0.1}^{[10]}$ & 1.68$\times10^{-7}$  & $\widetilde{\gamma}_{0.25}^{[10]}$   &   1.37$\times10^{-6}$  &  $\widetilde{\gamma}_{0.5}^{[10]}$  &  3.53$\times10^{-6}$ \\
  \hline
  
  $\widehat{\gamma}_{0.75}^{[10]}$   &   5.39$\times10^{-6}$ &  $\widehat{\gamma}_{0.9}^{[10]}$  &  1.34$\times10^{-5}$ & $\widehat{\gamma}_{0.95}^{[10]}$  &  2.51$\times10^{-5}$ & $\widehat{\gamma}_{0.975}^{[10]}$   &   4.57$\times10^{-5}$ \\
  
  $\widetilde{\gamma}_{0.75}^{[10]}$  &    8.46$\times10^{-6}$  &  $\widetilde{\gamma}_{0.9}^{[10]}$  &  2.05$\times10^{-5}$ & $\widetilde{\gamma}_{0.95}^{[10]}$   &  3.76$\times10^{-5}$  &  $\widetilde{\gamma}_{0.975}^{[10]}$   &   6.75$\times10^{-5}$ \\
    \hline

\end{tabular}
}}}
\label{tab:table2}
\end{table}

\begin{table}
\caption{$AMSE$ values when $F$ is weibull and $h = h^*$ or $h^{**}$}
{{\resizebox{\textwidth}{!}{\begin{tabular}{|c|c|c|c|c|c|c|c|}
  \hline
   &   $AMSE$  &  &  $AMSE$ &  &   $AMSE$  &  &  $AMSE$\\
  \hline
  $\widehat{w}_{0.05}^{[1/2]}$   &   4.11$\times10^{-2}$ & $\widehat{w}_{0.1}^{[1/2]}$  &  5.68$\times10^{-3}$ & $\widehat{w}_{0.25}^{[1/2]}$   &  3.85$\times10^{-4}$ & $\widehat{w}_{0.5}^{[1/2]}$   &  4.25$\times10^{-5}$ \\
  
   $\widetilde{w}_{0.05}^{[1/2]}$  &  4.20$\times10^{-2}$   & $\widetilde{w}_{0.1}^{[1/2]}$    &  5.93$\times10^{-3}$  &  $\widetilde{w}_{0.25}^{[1/2]}$    &  4.30$\times10^{-4}$  &  $\widetilde{w}_{0.5}^{[1/2]}$  & 5.50$\times10^{-5}$ \\
     \hline
     
  $\widehat{w}_{0.75}^{[1/2]}$   &    9.22$\times10^{-6}$  & $\widehat{w}_{0.9}^{[1/2]}$   &  3.31$\times10^{-6}$ & $\widehat{w}_{0.95}^{[1/2]}$   &   3.59$\times10^{-7}$ &  $\widehat{w}_{0.975}^{[1/2]}$  &  2.67$\times10^{-6}$  \\
  
  $\widetilde{w}_{0.75}^{[1/2]}$   &  1.53$\times10^{-5}$ &   $\widetilde{w}_{0.9}^{[1/2]}$  &  8.61$\times10^{-6}$   &  $\widetilde{w}_{0.95}^{[1/2]}$   &   7.68$\times10^{-6}$ & $\widetilde{w}_{0.975}^{[1/2]}$   &  7.90$\times10^{-6}$ \\
  \hline  \hline
  
  $\widehat{w}_{0.05}^{[10]}$   &   1.20$\times10^{-4}$ & $\widehat{w}_{0.1}^{[10]}$   &   2.65$\times10^{-4}$  & $\widehat{w}_{0.25}^{[10]}$   &   9.00$\times10^{-4}$ &  $\widehat{w}_{0.5}^{[10]}$   &  3.40$\times10^{-3}$ \\
  
  $\widetilde{w}_{0.05}^{[10]}$    &  1.11$\times10^{-4}$  &  $\widetilde{w}_{0.1}^{[10]}$    &  2.23$\times10^{-4}$ &  $\widetilde{w}_{0.25}^{[10]}$   & 4.79$\times10^{-4}$  &  $\widetilde{w}_{0.5}^{[10]}$ &  2.34$\times10^{-3}$ \\
  \hline
  
  $\widehat{w}_{0.75}^{[10]}$   &   1.38$\times10^{-2}$ &  $\widehat{w}_{0.9}^{[10]}$  &  5.49$\times10^{-2}$ & $\widehat{w}_{0.95}^{[10]}$    &  0.135 &  $\widehat{w}_{0.975}^{[10]}$  &  0.309 \\
  
  $\widetilde{w}_{0.75}^{[10]}$  &  1.33$\times10^{-2}$  &  $\widetilde{w}_{0.9}^{[10]}$  &  5.96$\times10^{-2}$   & $\widetilde{w}_{0.95}^{[10]}$   &   0.153 &  $\widetilde{w}_{0.975}^{[10]}$   &   0.360 \\
    \hline
\end{tabular}}}}
\label{tab:table3}
\end{table}

\begin{table}
\caption{$AMSE$ values when $F$ is beta and $h = h^*$ or $h^{**}$}
{{\resizebox{\textwidth}{!}{\begin{tabular}{|c|c|c|c|c|c|c|c|}
  \hline
   &   $AMSE$  &  &  $AMSE$ &  &   $AMSE$  &  &  $AMSE$\\
  \hline
  $\widehat{\beta}_{0.05}^{[1/2,1/2]}$  &  7.61$\times10^{-3}$ & $\widehat{\beta}_{0.1}^{[1/2,1/2]}$ & 1.20$\times10^{-3}$  &  $\widehat{\beta}_{0.25}^{[1/2,1/2]}$   &  1.42$\times10^{-4}$ &$\widehat{\beta}_{0.5}^{[1/2,1/2]}$  &  1.46$\times10^{-4}$ \\
  
  $\widetilde{\beta}_{0.05}^{[1/2,1/2]}$   &   7.77$\times10^{-3}$ & $\widetilde{\beta}_{0.1}^{[1/2,1/2]}$   &   1.25$\times10^{-3}$  & $\widetilde{\beta}_{0.25}^{[1/2,1/2]}$   &   1.54$\times10^{-4}$ & $\widetilde{\beta}_{0.5}^{[1/2,1/2]}$  & 1.11$\times10^{-4}$ \\
  \hline
  
  $\widehat{\beta}_{0.75}^{[1/2,1/2]}$   &   2.35$\times10^{-3}$ & $\widehat{\beta}_{0.9}^{[1/2,1/2]}$   &   0.167  & $\widehat{\beta}_{0.95}^{[1/2,1/2]}$   &   4.57  &  $\widehat{\beta}_{0.975}^{[1/2,1/2]}$  &  127 \\
  
  $\widetilde{\beta}_{0.75}^{[1/2,1/2]}$  &  1.61$\times10^{-3}$  & $\widetilde{\beta}_{0.9}^{[1/2,1/2]}$   &   0.116  & $\widetilde{\beta}_{0.95}^{[1/2,1/2]}$   &   3.17 & $\widetilde{\beta}_{0.975}^{[1/2,1/2]}$  &  87.9 \\
  \hline  \hline

$\widehat{\beta}_{0.05}^{[2,5]}$  &    1.55$\times10^{-4}$ &$\widehat{\beta}_{0.1}^{[2,5]}$   &   2.70$\times10^{-4}$ & $\widehat{\beta}_{0.25}^{[2,5]}$   &   5.35$\times10^{-4}$ & $\widehat{\beta}_{0.5}^{[2,5]}$    &  1.15$\times10^{-3}$ \\

$\widetilde{\beta}_{0.05}^{[2,5]}$   &   1.83$\times10^{-4}$ &  $\widetilde{\beta}_{0.1}^{[2,5]}$    &  2.18$\times10^{-4}$ & $\widetilde{\beta}_{0.25}^{[2,5]}$  &   2.91$\times10^{-4}$  &$\widetilde{\beta}_{0.5}^{[2,5]}$ &  5.67$\times10^{-4}$\\
\hline

$\widehat{\beta}_{0.75}^{[2,5]}$   &   3.09$\times10^{-3}$  &  $\widehat{\beta}_{0.9}^{[2,5]}$ &  1.01$\times10^{-2}$  &  $\widehat{\beta}_{0.95}^{[2,5]}$   &   2.41$\times10^{-2}$ & $\widehat{\beta}_{0.975}^{[2,5]}$   &  5.70$\times10^{-2}$ \\

$\widetilde{\beta}_{0.75}^{[2,5]}$   &   1.84$\times10^{-3}$  &  $\widetilde{\beta}_{0.9}^{[2,5]}$   &  7.22$\times10^{-3}$ &  $\widetilde{\beta}_{0.95}^{[2,5]}$   &   1.88$\times10^{-2}$  &  $\widetilde{\beta}_{0.975}^{[2,5]}$  &  4.72$\times10^{-2}$ \\
\hline\hline

 $\widehat{\beta}_{0.05}^{[5,2]}$  &  7.18$\times10^{-5}$ & $\widehat{\beta}_{0.1}^{[5,2]}$   &   1.38$\times10^{-4}$  &  $\widehat{\beta}_{0.25}^{[5,2]}$  &  4.20$\times10^{-4}$ & $\widehat{\beta}_{0.5}^{[5,2]}$   &   1.72$\times10^{-3}$  \\

 $\widetilde{\beta}_{0.05}^{[5,2]}$    &  6.57$\times10^{-5}$  &  $\widetilde{\beta}_{0.1}^{[5,2]}$   &  1.16$\times10^{-4}$ &  $\widetilde{\beta}_{0.25}^{[5,2]}$  &    2.84$\times10^{-4}$  &  $\widetilde{\beta}_{0.5}^{[5,2]}$   &   6.73$\times10^{-4}$\\
  \hline
  
  $\widehat{\beta}_{0.75}^{[5,2]}$    &  9.66$\times10^{-3}$ & $\widehat{\beta}_{0.9}^{[5,2]}$  &  6.68$\times10^{-2}$ &  $\widehat{\beta}_{0.95}^{[5,2]}$  &  0.261 &$\widehat{\beta}_{0.975}^{[5,2]}$  & 0.981  \\
  
  $\widetilde{\beta}_{0.75}^{[5,2]}$   &   4.23$\times10^{-3}$  & $\widetilde{\beta}_{0.9}^{[5,2]}$  & 4.01$\times10^{-2}$  &  $\widetilde{\beta}_{0.95}^{[5,2]}$   &   0.171  &  $\widetilde{\beta}_{0.975}^{[5,2]}$  & 0.676  \\
  \hline
\end{tabular}}}}
\label{tab:table4}
\end{table}

\section{Bias reduction and simulation study}
\label{sec:Bias reduction and simulation study}

As discussed in Section 3, the direct estimator has good performance, especially small variance.  If we could reduce the bias, we can get a better estimator.  There are many bias reduction methods, but the variance reduction is quite difficult.  In this section, we discuss reduction of the bias of $\widehat{H}(x_0)$.  The bias term is complicated, but if we use $4$-$th$ order kernel, we have $A_{1,2}=0$.  Thus we can reduce the convergence order of the bias from $O(h^2)$ to $O(h^4)$.  A simple way to construct the $4$-$th$ order kernel is proposed by Jones \cite{jones1997comparison}.  However the estimator takes negative value in some cases, though the hazard ratio is always non-negative. On the other hand, Terrell and Scott \cite{terrell1980improving}'s method reduces asymptotic bias without loss of non-negativity.  We propose the following modified direct hazard estimator
$$\widehat{H}^I(x) =\{\widehat{H}_{h}(x)\}^{\frac{4}{3}} \{\widehat{H}_{2h}(x)\}^{-\frac{1}{3}},$$
where $\widehat{H}_{h}(\cdot)$ and $\widehat{H}_{2h}(\cdot)$ are the direct estimator with the bandwidth $h$ and $2h$, respectively.  By \cite{terrell1980improving} we can get an asymptotic bias of $\widehat{H}^I(x_0)$
$$\frac{2 a_2^2(x_0) -4 a_4(x_0) H(x_0)}{H(x_0)} h^4,$$
where $a_2(x_0)$ and $a_4(x_0)$ are given in Theorem 2.3.  Whereas, the order of the variance is
$$Var[\widehat{H}^I(x_0)] =Var\left[\frac{4}{3} \widehat{H}_{h}(x_0) -\frac{1}{3}\widehat{H}_{2h}(x_0) \right] +O(n^{-1}) = O((n h)^{-1}),$$
and then the order is same as the original.  The order of the optimal bandwidth is $O(h^{-1/9})$, and the optimal convergence rate is improved to $O(n^{-{8}/{9}})$.  We confirm this in numerical simulations.

Table 5 shows the results of 10,000 times simulation of improvement rate (\%) of 
$$1-\{MISE[\widehat{H}^I]/MISE[\widehat{H}]\}),$$
when $F(\cdot)$ is the beta $B(r,s)$ and a sample size is $n=100$.  Their kernel functions are the Epanechnikov. The bandwidth $h^{\dagger}$ of $\widehat{H}$ is the asymptotically optimal in the sense of $MISE[\widehat{H}]$ which is available from Theorem 2.1 and that of $\widehat{H}^{I}$ is ${h^{\dagger}}^{5/9}$ modified to be  optimal in the sense of the convergence rate. Table 5 shows that how well the bias reduction works, and $\widehat{H}^I(x)$ is superior.  To obtain an explicit form of the asymptotic variance $Var[\widehat{H}^I(x)]$ is quite difficult and so we postpone it in future work.

\begin{table}
\caption{Improvement rate (\%) of $MISE[\widehat{H}^I]$ - $F$ is $B(r,s)*\sigma$ , $h = h^*$ or $h^{**}$}
{{\resizebox{\textwidth}{!}{\begin{tabular}{||c||c|c|c|c|c||}
  \hline
  Improvement rate (\%) & ~~~$\sigma=4^{-1}$~~~ & ~~~$\sigma=1$~~~ & $\sigma=4$ & ~~~$\sigma=4^2$~~~ & ~~~$\sigma=4^3$~~~\\
  \hline\hline
  $r=1/2, s=1/2$ & 0.301 & 1.25 & 1.52 & 1.64 & 1.78 \\
  \hline
  $r=1, s=1$ & 0.593 & 1.05 & 1.24 & 1.38 & 1.30 \\
  \hline
  $r=3, s=3$ & 6.74 & 8.35 & 8.27 & 8.28 & 8.35 \\
  \hline
  $r=2, s=5$ & 12.9 & 13.9 & 15.0 & 14.8 & 15.0 \\
  \hline
  $r=5, s=2$ & 3.96 & 5.10 & 5.11 & 5.05 & 5.05 \\
  \hline
\end{tabular}}}}
\label{tab:table5}
\end{table}

\section{Appendices: Proofs of Theorems}\label{appendices}

\noindent\textbf{Proof of Theorem 2.1}\medskip 

For the sake of simplicity, we use the following notations
\begin{eqnarray*}
&&\eta(z) = \int_{-\infty}^z F(u) du ,\hspace{10mm}\gamma(z) = \int_{-\infty}^z \eta(u) du, \\
&&M(z) = z -\eta(z) ,\hspace{10mm}m(z) =M'(z) =1-F(z).
\end{eqnarray*}
To begin with, we consider the following stochastic expansion of the direct estimator
\begin{eqnarray*}
&&\widehat{H}(x_0) \\
&=&\frac{1}{h} \int K\left( \frac{M(y) -M(x_0) }{h} \right) dF_n(y) \\
&&+ \frac{1}{h^2} \int K'\left( \frac{M(y) -M(x_0) }{h} \right) \{[\eta(y) -\eta_n(y)] -[\eta(x_0) -\eta_n(x_0)] \} dF_n(y) \\
&&+ \frac{1}{h^3} \int K''\left( \frac{M(y) -M(x_0) }{h} \right) \{[\eta(y) -\eta_n(y)] -[\eta(x_0) -\eta_n(x_0)] \}^2 dF_n(y) \\
&&+ \cdots \\
&=& J_1 + J_2 + J_3 + \cdots \hspace{70mm}({\rm say}).
\end{eqnarray*}

Since $J_1$ is a sum of $i.i.d.$ random variables, the expectation is obtained directly
\begin{eqnarray*}
E[J_1]&=& E \left[ \frac{1}{h} \int K\left( \frac{M(y) -M(x_0) }{h} \right) dF_n(y) \right] \\
&=& \frac{1}{h} \int K\left( \frac{M(y) -M(x_0) }{h} \right) f(y) dy \\
&=& \int K(u) \left[\frac{f}{1-F}\right] (M^{-1} (M(x_0) +h u)) du \\
&=& \left[\frac{f}{1-F}\right] (x_0) +\frac{h^2}{2} A_{1,2} \left[\frac{(1 -F)\{(1-F)f'' +4f f' \} +3f^3}{(1 -F)^5}\right](x_0) + O(h^4).
\end{eqnarray*}
Combining the following second moment
\begin{eqnarray*}
&& \frac{1}{n h^2} \int K^2\left( \frac{M(y) - M(x_0) }{h} \right) f(y) dy \\
&=& \frac{1}{n h}\int K^2(u) \left[\frac{f}{1-F}\right] (M^{-1} (M(x_0) +h u)) du \\
&=& \frac{1}{n h}\frac{f}{1-F} (x_0) A_{2,0} +O(n^{-1}),
\end{eqnarray*}
we have the variance as follows
$$
V[J_1] = \frac{1}{n h}\left[\frac{f}{1-F}\right] (x_0) A_{2,0} +O(n^{-1}).
$$

Next, we consider the following representation of $J_2$
$$J_2 = \frac{1}{n^2 h^2} \sum_{i=1}^n \sum_{j=1}^n K'\left( \frac{M(X_i) -M(x_0) }{h} \right) Q(X_i,X_j),$$
where 
$$Q(x_i,x_j) =[\eta(x_i) -J(x_i -x_j)] -[\eta(x_0) -J(x_0 -x_j)]$$
and
$$J (y-X_i) = I(X_i \le y) (y- X_i).$$
Using the conditional expectation, we get the following equation
\begin{eqnarray*}
E[J_2] &=& \frac{1}{n h^2} \sum_{j=1}^n E\left[ K'\left( \frac{M(X_i) -M(x_0) }{h} \right) Q(X_i,X_j) \right] \\
&=& \frac{1}{n h^2} E\left[ K'\left( \frac{M(X_i) -M(x_0) }{h} \right) E\left[ \sum_{j=1}^n Q(X_i,X_j) \Bigm| X_i \right] \right] \\
&=&  \frac{1}{n h^2} E\left[ K'\left( \frac{M(X_i) -M(x_0) }{h} \right) \left\{ \eta(X_i) -\left[ \eta(x_0) -J(x_0-X_i) \right] \right\} \right] \\
&=& \frac{1}{n h} \int K'\left( u \right) \left\{\eta(M^{-1} (M(x_0) +h u)) -\eta(x_0) +J(x_0-M^{-1} (M(x_0) +h u))  \right\} \\
&&\hspace{20mm}\times\frac{f}{1-F} (M^{-1} (M(x_0) +h u)) du\\
&=& \frac{1}{n h} \int K'\left( u \right) O(h u) \frac{f}{1-F} (x_0) du = O\left(\frac{1}{n}\right).
\end{eqnarray*}
Next we have
\begin{eqnarray*}
J_2^2 &=& \frac{1}{n^4 h^4} \sum_{i=1}^n \sum_{j=1}^n \sum_{k=1}^n \sum_{l=1}^n K'\left( \frac{M(X_i) -M(x_0) }{h} \right) K'\left( \frac{M(X_k) -M(x_0) }{h} \right) \\
&&\hspace{15mm}\times Q(X_i,X_j) Q(X_k,X_l) \\
&=& \frac{1}{n^4 h^4} \sum_{i=1}^n \sum_{j=1}^n \sum_{k=1}^n \sum_{l=1}^n J_2^2 (i,j,k,l)\hspace{30mm}({\rm say}).
\end{eqnarray*}
Taking a conditional expectation, we obtain that if all of $(i,j,k,l)$ are different, 
$$E[J_2^2 (i,j,k,l)] = E\left[ E\left\{J_2^2 (i,j,k,l) | X_i,X_k\right\} \right] = 0,$$
and
\begin{eqnarray*}
&&E[J_2^2 (i,j,k,l)]=0 \hspace{10mm}(\text{if $i=j$ and all of $(i,k,l)$ are different}),\\
&&E[J_2^2 (i,j,k,l)]=0  \hspace{10mm}(\text{if $i=k$ and all of $(i,j,l)$ are different}),\\
&&E[J_2^2 (i,j,k,l)]=0 \hspace{10mm}(\text{if $i=l$ and all of $(i,j,k)$ are different}).
\end{eqnarray*}
Then the terms in which `$j=l$ and all of $(i,j,k)$ are different' are the main term of $E[J_2^2]$.  If $j=l$ and all of $(i,j,k)$ are different, we have 
\begin{eqnarray*}
&&E[J_2^2 (i,j,k,l)] \\
&=& \frac{n(n-1)(n-2)}{n^4 h^4} E\biggl[ K'\left( \frac{M(X_i) -M(x_0) }{h} \right) K'\left( \frac{M(X_k) -M(x_0) }{h} \right)\\
&&\hspace*{35mm}\times Q(X_i,X_j) Q(X_k,X_j) \biggl].
\end{eqnarray*}
Using the conditional expectation of $Q(X_i,X_j) Q(X_k,X_j)$ given $X_i$ and $X_k$, we find
\begin{eqnarray*}
&&E\Bigl[ E\Bigl\{ Q(X_i,X_j) Q(X_k,X_j) \Bigm| X_i, X_k\Bigl\} \Bigl] \\
&=&E\Bigl[ \eta(X_i)\eta(x_0) +\eta(X_k)\eta(x_0) -\eta^2(x_0) +2\gamma(x_0) -\eta(X_i)\eta(X_k) \\
&&\hspace*{10mm} -(x +X_i -2\min(x ,X_i))\eta(\min(x ,X_i)) -2\gamma(\min(x ,X_i)) \\
&&\hspace*{10mm} -(x +X_k -2\min(x ,X_k))\eta(\min(x ,X_k)) -2\gamma(\min(x ,X_k)) \\
&&\hspace*{10mm} +(X_i +X_k -2\min(X_i ,X_k))\eta(\min(X_i ,X_k)) +2\gamma(\min(X_i ,X_k)) \Bigl].
\end{eqnarray*}
Therefore, the entire expectation of the last row is
\begin{eqnarray*}
&&E\biggl[ K'\left( \frac{M(X_i) -M(x_0) }{h} \right) K'\left( \frac{M(X_k) -M(x_0) }{h} \right) \\
&&\hspace*{10mm}\times (X_i +X_k -2\min(X_i ,X_k))\eta(\min(X_i ,X_k)) +2\gamma(\min(X_i ,X_k)) \Bigl]\\
&=&\int \Biggl[ \int_{\infty}^y K'\left( \frac{M(w) -M(x_0) }{h} \right) K'\left( \frac{M(y) -M(x_0) }{h} \right)\\
&&\hspace*{40mm}\times \left\{ (-w+y)\eta(w) + 2\gamma(w) \right\} f(y) dw \\
&&+\int_{y}^{\infty} K'\left( \frac{M(w) -M(x_0) }{h} \right) K'\left( \frac{M(y) -M(x_0) }{h}\right)\\
&&\hspace*{40mm}\times\left\{ (w-y)\eta(y) + 2\gamma(y) \right\} f(w) dw \Biggr] f(y) dy. \\
\end{eqnarray*}
Finally we get
\begin{eqnarray*}
&&E\biggl[ K'\left( \frac{M(X_i) -M(x_0) }{h} \right) K'\left( \frac{M(X_k) -M(x_0) }{h} \right) \\
&&\hspace{10mm}\times \{X_i +X_k -2\min(X_i ,X_k)\}\eta(\min(X_i ,X_k)) +2\gamma(\min(X_i ,X_k)) \Bigl]\\
&=& -h^2 \int K'\left( \frac{M(y) -M(x_0) }{h} \right) f(y) dy\\
&&\hspace{10mm}\times \Biggl( \Biggm[W\left( \frac{M(y) -M(x_0) }{h}\right)  \biggl( \left\{ \eta(x_0) +(-x_0 +y)F(x_0) \right\} \left[\frac{f}{m^2}\right](x_0) \\
&&\hspace{10mm} +\{(-x_0 +y) \eta(x_0) +2\gamma(x_0)\} \left[\frac{f'm -fm'}{m^3}\right](x_0) \biggl) +O(h)\Biggm] \\
&&\hspace{10mm}  +\Biggm[\left(1 -W\left( \frac{M(y) -M(x_0) }{h} \right) \right) \\
&&\hspace{10mm} \times\left( \eta(y) \frac{f}{m^2}(x_0) +\{(x_0-y) \eta(y)+ 2\gamma(y)\}\left[ \frac{f'm -fm'}{m^3}\right](x_0) \right) +O(h) \Biggm] \Biggl)\\
&=& h^4 \left[ \frac{Ff^2}{m^4}(x_0) +\left\{\frac{f'm -fm'}{m^3} \left( 2\eta \frac{f}{m^2} +2\gamma \frac{f'm -fm'}{m^3} \right)\right](x_0) \right\} +O(h^5).
\end{eqnarray*}
After similar calculations of the other terms, if $j=l$ and all of $(i,j,k)$ are different, we get
$$E[J_2^2 (i,j,k,l)] = O\left( \frac{1}{n}\right).$$
In addition, it is easy to see that the expectations of the other combinations of $(i,j,k,l)$ are $o(n^{-1})$. As a result, we have
$$
E[J_2^2] = O(n^{-1}) \hspace{3mm}\text{and}\hspace{3mm} V[J_2] = O(n^{-1}).
$$

The moments of $J_3$ is also obtained in a similar manner and we find that $J_3+ \cdots$ is negligible. The main bias of $\widehat{H}(x_0)$ comes from $J_1$. By Cauchy-Schwarz inequality, we find that the main term of the variance comes from $J_1$.  Now, we get the $AMSE$ of the direct estimator and the proof of Theorem 2.1.\\

\noindent\textbf{Proof of Theorem 2.2 and 2.3}\medskip

It follows from the previous discussion that
\begin{eqnarray*}
&&\sqrt{n h} \left\{ \widehat{H}(x_0) -\left[\frac{f}{1-F}\right](x_0) \right\} \\
&=& \sqrt{n h} \left\{ J_1 -\left[\frac{f}{1-F}\right](x_0) \right\} +o_P(1) \\
&=& B_1 + \sqrt{n h} \left\{ J_1 -\left[\frac{f}{1-F}\right](x_0)-h^2 B_1 \right\} +o_P(1).
\end{eqnarray*}
Since $J_1$ is a sum of $i.i.d$. random variables and the expectation of the second term is $o(1)$, the asymptotic normality of Theorem 2.2 holds. Furthermore, these arguments show that
\begin{eqnarray*}
E \left[ \widehat{H}(x_0) \right] &=& E[J_1] +O(n^{-1}) \\
&=& \int K(u) \left[\frac{f}{1-F}\right] (M^{-1} (M(x_0) +h u)) du +O(n^{-1}),
\end{eqnarray*}
and we get directly Theorem 2.3 by the Taylor expansion.

\end{document}